\documentclass[twocolumn]{autart}    

\usepackage{natbib}
\setcitestyle{authoryear,open={(},close={)}}

\usepackage{url}
\usepackage{graphicx}          
\usepackage{amsmath,amssymb,setspace,calrsfs}

\newif\iffinalfonts
\iffinalfonts
\usepackage{mbboard}
\else
\newcommand{\bbPi}{\Pi}
\newcommand{\bbSigma}{\Sigma}
\fi

\usepackage{pgfplotstable,pstool,psfrag}
\usepackage{float}
\usepackage{tabularx}
\usepackage{ltablex}
\usepackage{booktabs}
\usepackage{picins}

\newcommand{\R}{\mathbb{R}}
\newcommand{\N}{\mathbb{N}}

\newcommand{\probability}{\mathbb{P}}

\def\defeq{\mathrel{\mathop:}=}

\newtheorem{prp}[thm]{Proposition}

\begin{document}

\begin{frontmatter}

\title{On the Ergodic Control of Ensembles\thanksref{footnoteinfo}}
\thanks[footnoteinfo]{This paper was not presented at any IFAC
meeting.
}

\author[UNICAMP]{Andr{\' e} R.~Fioravanti}\ead{fioravanti@fem.unicamp.br},
\author[IBM]{Jakub Mare\v{c}ek}\ead{jakub.marecek@ie.ibm.com},
\author[IBM,UCD,IMPERIAL]{Robert N. Shorten}\ead{robert.shorten@ucd.ie},
\author[UNICAMP]{Matheus Souza}\ead{msouza@fee.unicamp.br},
\author[PASSAU]{Fabian R. Wirth}\ead{fabian.wirth@uni-passau.de}

\address[UNICAMP]{University of Campinas, Brazil}  
\address[IBM]{IBM Research -- Ireland, in Dublin, Ireland}             
\address[UCD]{University College Dublin, Dublin, Ireland}        
\address[IMPERIAL]{Imperial College London, South Kensington, UK}        
\address[PASSAU]{University of Passau, Germany}        


\begin{abstract}
Across smart-grid and smart-city application domains,
there are many problems where an ensemble
of agents is to be controlled such that both the
aggregate behaviour and individual-level perception
of the system's performance are acceptable.
In many applications, traditional PI control is used
to regulate aggregate ensemble performance.
Our principal contribution in this note is to demonstrate that
PI control may not be always suitable for this purpose, and in
some situations may 
lead to a loss of ergodicity for
closed-loop systems. Building on this observation, a theoretical framework is proposed
to both analyse and design control systems for the regulation
of large scale ensembles of agents with a {\em
probabilistic intent}. Examples are given to illustrate
our results. \end{abstract}

\end{frontmatter}

\maketitle

\section{Introduction}
At a very high level, smart-city related research concerns designing systems that endeavour to make the best use of limited resources across a number of
domains (energy, transport, water etc).  While classical control has much to offer in such application areas,
there are aspects and peculiarities of many of these applications that require
practitioners in Control Theory to explore new types of theoretical and practical challenges.
Roughly speaking, classical control is typically concerned with regulating a single system, such that the system achieves a desired behaviour in an optimal way.
In contrast, in many smart city applications, we are interested in allocating a resource among a population of
agents. These might be humans or algorithmic processes that bid for access to a resource in some probabilistic
manner (for example, access to part of a road network). In such applications, both the experience of the individual, and
their aggregate effect are important. Furthermore, in smart-city applications, we typically wish to control and
influence the behaviour of large-scale populations,  where the number of agents varies over time and
is not known with certainty. Additionally, there are limits to the observability of
such systems and data sets are often obtained in a closed-loop fashion; that is, operator's decisions are often reflected in the data sets.
Finally, a fundamental difference between classical control and smart-grid and smart-city control
is the need to study the effect of a common control signal on the individual agent and its long term access to a constrained
resource.  Among all of these fundamental differences, it is this last issue
that is perhaps most alien to the classical control theorist, and yet the issue that is perhaps the most pressing in real-life applications,
since the need for predictability, at the level of individual agents,
underpins an operator's ability to write economic contracts.

In this paper,
our starting point is the observation that many problems that are considered in smart-grids and smart-cities
 can be cast in a framework where a large number of agents, such as people, cars, or machines, often with unknown objectives,
 compete for a limited resource. The challenge of allocating this resource in a manner that is not wasteful,
 which gives an optimal return on the use of the resource for society,
 and which, in addition, gives a guaranteed level of service to each of the agents competing for that resource,
 gives rise to a whole host of problems, which in principle are best addressed in a control-theoretic manner.
%
From the perspective of a control engineer, this statement
can be decomposed into three objectives; two of which are familiar in control, and
the other one constitutes a relatively new consideration.
Our first objective is to fully utilise the resource, which is a regulation problem.
Second, we would then like to make optimal use of the resource.
While both of these objectives are concerned with the aggregate
behaviour of an agent population, they make no attempt to control the manner in which the agents orchestrate their behaviour
to achieve this aggregate effect.
Our third objective thus focuses on the effects of the control on the microscopic properties of the agent population.
Ultimately, this third objective can be phrased in terms of properties of the stochastic process capturing the share of the resource
that is allocated to an {\em individual agent}. For example, we may wish that each agent, on average, receives a fair share of the
resource over time, or, at a much more fundamental level, we wish the average allocation of the resource to each agent over time to
be a stable quantity that is entirely predictable and
which does not depend on initial conditions, and which is not sensitive to noise entering the system. From the point of view of the
design of the feedback system, these latter concerns are related to the existence of the {\em unique invariant measure} that governs
the distribution of the resource amongst the agents in the long run.
Thus, the design of feedback systems for deployment in multi-agent applications must consider not only the traditional notions of regulation and optimisation,
but also the guarantees concerning
the existence of this unique invariant measure. As we shall see, this is not a trivial task and many familiar control strategies,
in very simple situations,
do not necessarily give rise to feedback systems which possess all three of these features.

\subsection{Brief overview of related work}

Due to the intrinsically multidisciplinary nature of the problem studied in this paper,
related topics are discussed in several communities.
As stated before, several problems in smart cities and smart grid can be seen as the task of efficiently allocating some
resource among a population of agents, even though there is no formal and widely accepted definition for ``smart city'' \citep{IOT2014}.
Intelligent transportation systems (ITS), for instance, are one of the
many smart-city application areas that thrived in the last few years,
presenting rapid technological developments due to (i) the integration of transportation systems with the Internet and
(ii) the pressure for green -- environmentally friendly -- transport solutions.
Typical applications in ITS vary from monitoring and controlling transportation flows \citep{ITS2010} and $CO_2$ emissions \citep{Arieh13}
to assigning parking spaces \citep{Arieh14} and utilising vehicle-to-grid (V2G) energy supply \citep{EV2018,5560848}.
In the context of smart grids, problems related to real-time electricity pricing \citep{SG2010a},
real-time demand response \citep{SG2010,Hiskens2011}, and real-time interruptible loads management \citep{192993,Alagoz2013,Salsbury2013}
and V2G-related issues,
have also attracted interest over the last few years.
For these -- and many other -- applications, we present both conditions ensuring and ruling out ergodicity
in a certain closed-loop sense. Related work also appears in the context of multi-agent systems \citep{Agents07}.
Multi-agent systems arise, for instance, in the application areas stated above, and many references in the literature address
these problems using ideas based on {\em consensus} or {\em agreement} \citep{Blondel05,Nedic09}.
The consensus approach is very useful in practice, due to its strong links to utility maximization and fairness.
The concept of fairness can be related to {\em ergodicity}, 
as an ergodic dynamic behaviour implies several important properties that are necessary for fairness \citep{Mezic11}.
Finally, we note that our positive results in this paper are based on {\em iterated function systems}
\citep{elton1987ergodic,BarnsleyDemkoEltonEtAl1988,barnsley1989recurrent}. Iterated function systems (IFS), albeit
not well known in the control community, are a convenient and rich class of Markov processes.
As it will be seen in the sequel, a class of stochastic systems arising from the  dynamics of multi-agent interactions
can be modelled and analysed using IFS in a particularly natural way.
The use of IFS makes it possible to obtain strong stability guarantees for such stochastic systems.
To prove certain contrapositive statements, suggesting when such guarantees are impossible to provide,
we use novel coupling arguments of Hairer et al. \citep{hairer2011asymptotic}.
Our use of coupling arguments to prove such negative results is one of the first within Control Theory,
as far as we know.

\subsection{Paper Organisation and Contributions}

The main purpose of the present paper is to introduce and analyse the main features of a problem class which
is of great interest for smart cities applications. Our main contributions, apart from formulating the problem class itself, arise
from the study of the ergodic properties associated with output-regulation problems in such multi-agent settings. At a high level, we
show the following:
\begin{itemize}
\item
In regulating the behaviour of ensembles of agents,
feedback control with integral action may destroy the ergodic properties
of the closed-loop system, even when the agent behaviour is benign and despite the fact that regulation is achieved.
The importance of this result stems from the fact that ergodic behaviour is a fundamental property that is essential
for underpinning economic contracts, and for guaranteeing properties such
as fairness.
The result is hence important from a practical context, and is not merely of theoretical interest.
\item 
We present specific examples to illustrate the loss of ergodicity in very benign examples.
\item We show that for certain population types and filters, stable control action always results in ergodic behaviour.
In particular, we show this for linear and non-linear systems,
with both real-valued actions, and for actions constrained to a finite set.
\item A final minor contribution is to illustrate the use of results from the study of iterated function systems in designing controllers for certain classes of
dynamic systems.
\end{itemize}
Given the widespread interest in applications where the behaviour of ensembles of rational agents are orchestrated with
probabilistic intent,
we believe that these results may be of interest in these
and related areas.


The paper is organised as follows: In
  Section~\ref{sec:preliminaries}, we present our model. We
  then formulate the necessary concepts from the theory of Markov chains we
  will use, and recall the concept of
  coupling of invariant measures in order to state a necessary condition
  for ergodicity.
  In Section~\ref{sec:comments-pi-negative}, we present a negative result, which shows
  that ergodicity may fail whenever a standard PI controller is used in the loop.
In particular, the amount of a resource that can be used by a particular agent depends on the initial state of the controller.
In Section~\ref{sec:comments-pi-positive}, a positive result is
  obtained for stable linear controllers.
In Sections \ref{sec:nonlinear}, and \ref{sec:discrete},
we extend the results to non-linear systems and cases
where actions are constrained to a finite set,
respectively.

\begin{rem} A preliminary version of this paper has appeared in \citep{8263852}. The present paper extends
 beyond this preliminary version in several ways.
 First, full proofs are given.
 Second, additional positive results are developed for both linear and non-linear systems with both continuous and finite sets of actions.
 Finally, examples illustrating loss of ergodicity are presented.
\end{rem}

\section{Preliminaries}\label{sec:preliminaries}

We now develop the general setting of this paper. The objective here is to set out our modelling framework,
and to present basic results that can be useful in studying the properties of control strategies for ensembles.

\subsection{Notation}

In general, our notation employs the following rules:
Upper-case letters are used for matrices,
in caligraphic they denote groups, maps, operators, and in blackboard-bold
the probability operator $\mathbb{P}$ as well as sets, and spaces
while lower-case letters are used for vectors, scalars, and functions.
Subscripts are used to distinguish symbols; time-indexed symbols are followed by the time index in parentheses, as in $x(k)$.
Superscript is used for exponentiation and, in $(\cdot)^T$, for the transpose operation.
For a given measurable space $(\mathbb{X}, {\mathcal F})$, with
$\sigma$-algebra ${\mathcal F}$, $M(\mathbb{X})$ indicates the set of all probability measures over $\mathbb{X}$ and $\mathbb{X}^\infty$
denotes its associated path space, which consists of infinite right-sided sequences over $\mathbb{X}$: $x(0) \in \mathbb{X}, x(1) \in \mathbb{X}, x(2) \in \mathbb{X}, \ldots$.
Details of the notation are introduced locally in the sections, where it is first utilised.
For a complete listing of the symbols utilised, please see the supplementary material\footnote{\url{https://arxiv.org/abs/1807.03256}} on-line.

\subsection{Models}
	\label{sec:models}

We consider the problem of repeatedly distributing a limited resource among multiple agents, based on some information concerning the resource, which is provided by a central authority.
Throughout, we consider systems subject to several constraints.
First, the central authority does not observe the consumption of individual agents, but rather the total utilisation of the resource, or a filtered version thereof.
Based on the filtered measurements of the utilisation, the central authority
 provides information to the agents, sets the price of utilising the resource, or similar.
Second, the agents respond to information broadcast by the central authority, but have only limited communication capability otherwise.
Specifically, we assume no inter-agent communication.
Third, the agents have their own, private objectives.
That is, although they receive information from the central authority, they need not pick an action the authority would deem most appropriate. As we shall see, it will be convenient to encode the selfish response of an agent to the information in a probabilistic manner.
Finally, in our model, the agents may be limited to a choice from a finite set of possible requests for the resource. In an extreme case, the agents only have the possibility  to turn their utilisation on or off, i.e. $x_i\in \{ 0,1 \}$. In a more general setting, a
subpopulation might be able to choose their consumption from a continuous interval or via some local control.

\begin{figure}
  \centering
  \includegraphics[width=\columnwidth]{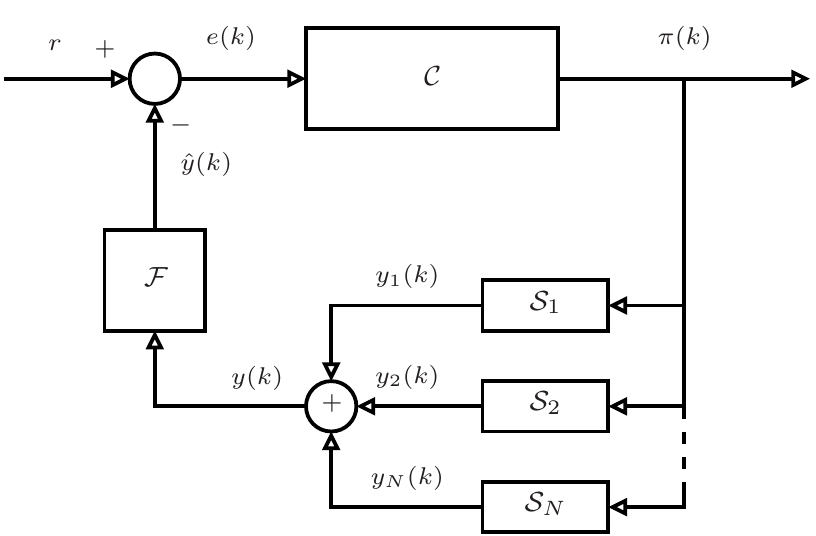}
  \caption{Feedback model.}\label{system}
\end{figure}

With these constraints in mind, we are interested in the closed-loop,
discrete time system depicted in Figure \ref{system},
comprising a controller, a number $N \in \N$ of agents, and a filter.
A controller  $\mathcal{C}$, which represents the central authority,
produces a signal $\pi(k)$ at time $k$.
In response, the agents,
modelled by systems $\mathcal{S}_1$, $\mathcal{S}_2$, \ldots, $\mathcal{S}_N$,
amend their use of the resource.
The internal state of the agents is
denoted by $x_i$, $i=1,\ldots,N$. In particular, $x_i(k)$ is a random variable.
We model the use of the resource $y_i(k)$ of agent $i$ at time $k$ as the
output of system ${\mathcal S}_i$, which is again a random variable. In the
remainder we will assume $y_i(k)$ and $\pi(k)$ are scalars, but generalisations are easy to obtain.
 The randomness can be a result of the inherent randomness in the
reaction of user $i$ to the control signal $\pi(k)$, or the response to a
control signal that is intentionally randomized \citep{Arieh13,Arieh14,marevcek2015signaling}.
The aggregate resource utilisation $y(k) \defeq \sum_{i = 1}^N y_i(k)$   at time $k$
is then also a random variable.
The controller may not have access to either $x_i(k)$, $y_i(k)$ or $y(k)$,
but only to
the error signal $e(k)$, which is the difference of
$\hat y(k)$, the output of a filter $\mathcal{F}$,
and $r$, the desired value of $y(k)$.
Further, we assume that the controller has its private state $x_c(k) \in \R^{n_c}$.
The controller aims to regulate the system by providing a signal $\pi(k) \in \bbPi \subseteq
\R$ at time $k$; here, $\bbPi$ denotes the set of admissible broadcast control signals.
In the simpler static case,
the signal $\pi(k)$ is a function of an error
signal $e(k)$ and the controller state $x_c(k)$,
whose range is $\bbPi$.

In the case of purely discrete agents, the non-de\-ter\-min\-is\-tic agent-specific response to the feedback signal $\pi(k)\in \bbPi$ can be modelled
by agent-specific and signal-specific probability distributions over certain agent-specific set of actions
$\mathbb{A}_i =\{
    a_1,\ldots,a_L \}\subset \R^{n_i}$,
where $\R^{n_i}$ can be seen as the space of
agent's $i$ private state $x_i$.
Assume that the set of possible resource demands of agent $i$ is
$\mathbb{D}_i$, where in the case that $\mathbb{D}_i$ is finite we denote
\begin{equation}
    \label{eq:5}
    \mathbb{D}_i := \{ d_{i,1},
    d_{i,2}, \ldots, d_{i,m_i}\}.
\end{equation}
In the general case, we assume there are $w_i \in \N$ state transition maps ${\mathcal W}_{ij}: \R^{n_i} \to \R^{n_i}$, $j=1,\ldots,w_i$
for agent $i$
    and $h_i \in \N$ output maps ${\mathcal H}_{i\ell}: \R^{n_i} \to
    \mathbb{D}_i$, $\ell= 1,\ldots,h_i$ for each agent $i$.
The evolution of the states and
 the corresponding demands then satisfy:
\begin{align}
\label{general}
x_i(k+1) & \in  \{  {\mathcal W}_{ij}(x_i(k)) \;\vert\; j = 1, \ldots, w_i\}, \\
y_i(k) & \in  \{ {\mathcal H}_{i\ell}(x_i(k)) \;\vert\; \ell = 1, \ldots, h_i\},
\end{align}
where the choice of agent $i$'s response at time $k$ is governed by probability
functions $p_{ij} : \bbPi \to [0,1]$, $j=1,\ldots,w_i$, respectively
$p'_{i\ell} : \bbPi \to [0,1]$, $\ell=1,\ldots,h_i$. Specifically, for each agent $i$, we have
for all $k\in\N$ that
\begin{subequations} \label{eq:problaws}
\begin{align}
&\mathbb{P}\big( x_i(k+1) = {\mathcal W}_{ij}(x_i(k)) \big) = p_{ij}(\pi(k)),\\
&\mathbb{P}\big( y_i(k) = {\mathcal H}_{i\ell}(x_i(k)) \big) = p'_{i\ell}(\pi(k)).
    \label{eq:7}
\intertext{Additionally, for all $\pi \in \Pi$, $i=1,\ldots,N$ it holds that}
&    \sum_{j=1}^{w_i} p_{ij}(\pi) = \sum_{\ell=1}^{h_i} p'_{i\ell}(\pi) = 1.
\end{align}
\end{subequations}
The final equality comes from the fact $p_{ij}$, $p'_{i\ell}$ are
probability functions.
We assume that,
conditioned on $\{ x_i(k)  \}, \pi(k)$,
the random variables $\{ x_i(k+1) \mid  i = 1,\ldots,N \}$ are
stochastically independent. The outputs $y_i(k)$ each depend on $x_i(k)$ and the
signal $\pi(k)$ only.
As we shall see, this general framework allows for surprisingly sharp results.

Our specific aim is to distribute the resource such
that we achieve the following goals almost surely, i.e. with probability $1$:
\begin{enumerate}
  \item \textbf{feasibility}:
Given an upper bound $r>0$ for the utilization of the resource, we require
for all $k\in\N$
    \begin{equation}
        \label{eq:1}
        \sum_{i = 1}^N y_i(k) = y(k) \leq r.
    \end{equation}
    More generally, the resource could be time-varying; for the purposes
of this paper it will assumed to be a constant quantity.
  \item \textbf{predictability}: for each agent $i$ there exists a
    constant $\overline{r}_i$ such that
    \begin{equation}
        \label{eq:3}
        \lim_{k\to \infty} \frac{1}{k+1} \sum_{j=0}^k y_i(j) = \overline{r}_i,
    \end{equation}
where this latter limit is independent of initial conditions.
\end{enumerate}
Further optional requirements may include: \textbf{fairness}, which could be
formulated by saying that all the $\overline{r}_i$ coincide, and
\textbf{optimality} so that the vector $\overline{r} =
\begin{bmatrix}
    \overline{r}_1 & \ldots & \overline{r}_N
\end{bmatrix}$ is a local optimum of an underlying optimization problem.
In addition, it is
also of interest to achieve the goals after a transient phase, i.e. for
all $k\geq K$, where $K$ is a constant.

While all of these goals are important from a practical perspective,
the principal property of interest in this paper is the goal of {\em predictability},
since this latter issue defines the ability of service providers to write economic contracts.
In order to analyze this more formally, we consider an augmented state space $\mathbb{X} \subset \R^d$, which captures the state of the controller, the filter, and the agents.
Denote by $M(\mathbb{X})$ the space of probability measures over
$\mathbb{X}$ with the standard $\sigma$-algebra.
The behaviour of the overall system in response to the signal $\pi(k)$ can be modelled as $P: \mathbb{X} \times \bbPi \to M(\mathbb{X})$.
In order to reason about the evolution of the state, consider the space $\mathbb{X}^{\infty}$ of one-sided infinite sequences of states,
known as the path space.
Further, we introduce
the space of probability measures on path space $\mathbb{X}^{\infty}$ with
the product $\sigma$-algebra.
Notice that for a particular combination of a filter ${\mathcal F}$,
controller ${\mathcal C}$, and population of agents,
the feedback loop can be modelled by an 
operator $\mathcal P_k: M(\mathbb{X}) \to M(\mathbb{X})$ and the associated dynamical system $(\mathcal P_k)_{k \in \N}$.
Our paper, at some level, asks what properties of ${\mathcal C}$ and ${\mathcal F}$ make $(\mathcal P_k)_{k \in \N}$ predictable
and what properties of ${\mathcal C}$  render it impossible for it to be predictable.
Predictability is related to asymptotic convergence in probability in terms of measures over the path space, i.e., $M(\mathbb{X}^{\infty})$,
and will be conveniently characterised in the next section in
terms of the existence of a {\em unique invariant measure},
in the language of Markov chains.

\subsection{Markov Chains and Iterated Function Systems}

The set-up described in the prequel resembles closely that of an iterated function system \citep{elton1987ergodic,BarnsleyDemkoEltonEtAl1988,barnsley1989recurrent}.
Iterated function systems are a class of stochastic dynamic systems, for which strong stability and convergence
results exist. We now introduce some notation and mention some of the most important results.

To begin, let $\bbSigma$ be a closed subset of $\R^n$ with the usual Borel
$\sigma$-algebra $\mathbb{B}(\bbSigma)$. We call the elements of $\mathbb{B}$
events. A Markov chain on $\bbSigma$ is a sequence
of $\bbSigma$-valued random
vectors $\{ X(k)\}_{k\in\N}$ with the Markov property, which is the equality of
a probability of an event conditioned on past events
and probability of the same event conditioned on the current state, \emph{i.e.}, we always have
\begin{align}
    \probability \left( X(k+1) \in \mathbb{G} \mid X(j)=x_{j}, \,  j=0, 1, \dots, k \right) \nonumber \\
	= \probability\big(X(k+1)\in \mathbb{G} \mid X(k)=x_{k}\big), \nonumber
\end{align}
	where $\mathbb{G}$ is an event
and $k \in \N$. We assume the Markov chain is time-homogeneous. The transition operator $P$ of the Markov chain is defined
for $x \in \bbSigma$, $\mathbb{G} \in \mathbb{B}$ by
\begin{equation*}
    P(x,\mathbb{G}) := \probability(X(k+1) \in \mathbb{G} \; \vert \; X(k) = x).
\end{equation*}
 If the initial condition $X_0$ is distributed according to an initial distribution
$\lambda$, we denote by $\probability_\lambda$ the probability measure induced on
the path space, i.e., space of sequences with values in $\bbSigma$. Conditioned on an
initial distribution $\lambda$, the random variable $X(k)$ is distributed
according to the measure $\lambda_k$ which is determined inductively by
\begin{equation}
    \label{eq:measureiteration}
    \lambda_{k+1}(\mathbb{G})  := \int_{\bbSigma} P(x, \mathbb{G}) \, \lambda_k(d x),
\end{equation}
for $\mathbb{G} \in \mathbb{B}$. A measure $\mu$ on $\bbSigma$ is called
invariant with respect to the Markov process $\{ X(k) \}$ if it is a fixed
point for the iteration described by \eqref{eq:measureiteration},
\emph{i.e.}, if $\mu P = \mu$. An invariant probability measure $\mu$ is called
attractive, if for every probability measure $\nu$ the sequence $\{
\lambda_k \}$ defined by \eqref{eq:measureiteration} with initial
condition $\nu$ converges to $\mu$ in distribution.
The existence of attractive invariant measures is
intricately linked to ergodic properties of the system.

With this background, our general problem considered in this paper is modelled as a Markov chain on a state space
representing all the system components. Let $\mathbb{X}_S$
be the product space of the state spaces of all agents. The
spaces $\mathbb{X}_F,\mathbb{X}_C$ contain the possible internal states for filter and
central controller. Our system thus evolves on the state space $\mathbb{X} :=
\mathbb{X}_S \times \mathbb{X}_C \times \mathbb{X}_F$,
as we will see in the following sections.

\subsection{Invariant Measures and Ergodicity}
\label{sec:ergodicity}

In our positive results, we consider a class of Markov chains that
are known as {\em iterated function systems} (IFSs). In an iterated
function system, we are given a set of maps $\{ f_j : \bbSigma \to \bbSigma
\; \vert \; j \in \mathbb{J} \}$, where $\mathbb{J}$ is a (finite or countable) index set. Associated to
these maps, there are probability functions $p_j: \bbSigma \to [0,1]$,
where $\sum_{j\in \mathbb{J}} p_j(x) = 1$ for
all $x \in \bbSigma$.
The state $X(k+1)$ at time $k+1$ is then given by $f_j(X(k))$ with probability $p_j(X(k))$,
where $X(k)$ is the state at time $k$.
Such IFSs have been studied in fractal image compression \citep{barnsley2013fractals},
congestion management in computer networking \citep{corless2016aimd} and transportation \citep{Marecek2016},
and to some extent, control theory at large \citep{Branicky1994,Branicky1998,Marecek2017}.

Sufficient conditions for the existence of a unique attractive
invariant measure can be given in terms of
 ``average contractivity''. This key notion can be traced back to
\citep{elton1987ergodic,BarnsleyDemkoEltonEtAl1988,barnsley1989recurrent}:
\begin{thm}[\cite{BarnsleyDemkoEltonEtAl1988}]
\label{Barnsley}
  Let $\bbSigma \subset \R^n$ be closed.
  Consider an IFS with a finite index set $\mathbb{J}$,  and globally
  Lipschitz continuous maps
  $f_j:\bbSigma \to \bbSigma$, $j\in \mathbb{J}$.
  Assume that the probability functions $p_j$ are globally Lipschitz continuous and
  bounded below by $\eta > 0$.
  If there exists a $\delta>0$  such that for all $x, y \in \bbSigma, x \not = y$
  \begin{align*}
    \sum_{j\in \mathbb{J}} p_j(x) \log \left( \frac{ \|{f_j (x) - f_j (y)}\| }{ \|x-y\| } \right) < -\delta < 0,
  \end{align*}
  then there exists an attractive (and hence unique) invariant probability
  measure
  $\mu$ for the IFS.
\end{thm}

We can combine Theorem~\ref{Barnsley} with a theorem by Elton
\citep{elton1987ergodic}, to obtain that for all (deterministic) initial conditions $x_0 \in
\bbSigma$ and continuous $g:\bbSigma\to \R$, the limit
\begin{equation}
\label{eq:ergodicprop}
   \lim_{k \to\infty} \frac{1}{k+1} \sum_{\nu=0}^k  g(X(\nu)) = \mathbb{E}_\mu(g)
\end{equation}
exists almost surely ($\probability_{x_0}$)  and is
independent of $x_0 \in \bbSigma$. The limit is given by the expectation
with respect to the invariant measure $\mu$.
For more general theorems,
the reader is referred to
\citep{elton1987ergodic,BarnsleyDemkoEltonEtAl1988,barnsley1989recurrent,barnsley2013fractals,stenflo2001markov,szarek2003invariant,steinsaltz1999,walkden2007invariance,barany2015iterated}
and especially two recent surveys
\citep{iosifescu2009iterated,stenflo2012survey}.

\begin{rem}
	\label{rem:baseline}
    From the point of view of applications in smart cities, the existence of
	the limit \eqref{eq:ergodicprop} is a minimum requirement. We want to avoid
    situations where the average allocation of resources to agents depends
    on their initial conditions, on possible initial conditions of
    controllers and filters, etc. In addition, it is desirable to shape
    the expected value so that an overall optimum is obtained.
    In fact, if Theorem~\ref{Barnsley} and
      \citep{elton1987ergodic} are applicable, then predictability (cf. Section \ref{sec:models} above) holds with probability one and the limit is independent
      of initial conditions. It is relatively easy to ensure that the
      probabilities for all state transition maps are always positive by
      designing rules for each agent. For the
    question of feasibility (cf. Section \ref{sec:models} above) the shaping of the expectation is essential.
\end{rem}

There is a vast literature on invariant measures and ergodic properties of stochastic systems.
In place of a definition of ergodicity, we summarise some of the main results following \citep{Hairer2006}:

\begin{prp}
\label{prop:ergodic-car}
Consider an IFS on the state space $\bbSigma$ with
invariant probability measure  $\mu$.
The following are equivalent:
\begin{enumerate}
\item[E1] \label{ergodic} $\mu$ is ergodic.
\item[E2]
\label{C511} every $\mu$-invariant set $S \subset \bbSigma$ is of $\mu$-measure 0 or 1, where
a measurable set $S \subset \bbSigma$ is $\mu$-invariant if $P(x, S) = 1$ for $\mu$-almost every $x \in S$.
\item[E3] \label{extremal} $\mu$ cannot be decomposed as $\mu = t \mu_1 +
  (1 - t) \mu_2$ with $t \in (0,1)$ for two invariant probability measures $\mu_1, \mu_2$.
                                \end{enumerate}
Additionally, 
the following implies that $\mu$ is ergodic:
\begin{enumerate}
\item[F1] \label{unique} the Markov process with transition operator $P$ has a unique invariant measure.
\end{enumerate}
Additionally, if $\mu$ is ergodic, the following holds:
\begin{enumerate}
\item[C1] \label{Breiman} $\lim_{k \to \infty} \frac{1}{k+1} \sum_{\nu =
    0}^{k} g(x(k)) = \mathbb{E}_\pi(g)$ almost surely, for
  $\mu$-almost all initial conditions.
\item[C2] \label{Orthogonal} any other distinct ergodic invariant measure is singular to $\mu$.
\end{enumerate}
\end{prp}

\begin{pf}
E1 is equivalent to E2 by Corollary 5.11 of \citep{Hairer2006}.
E1 is equivalent to E3 by Theorem 5.7 of \citep{Hairer2006}.
F1 implies E1 by Corollary 5.12 of \citep{Hairer2006}.
E1 implies C1 by Corollary 5.3  of \citep{Hairer2006}.
E1 implies C2 by Theorem 5.7 of \citep{Hairer2006}.
\end{pf}

We note that it is possible that multiple ergodic invariant measures exist
for a given system or Markov process. This case is however not of interest
for us. We call an IFS uniquely ergodic if it has an attractive, unique,
ergodic, invariant probability measure, such that
\eqref{eq:ergodicprop} holds for all deterministic initial
conditions. Note that Theorem~\ref{Barnsley} together with
\citep{elton1987ergodic} provide sufficient conditions for unique
ergodicity in this sense.

As we have discussed before this notion precisely guarantees the property
of predictability and can be used for an analysis of fairness and optimality.

\subsection{Ergodic Invariant Measures and Coupling}

In our negative results, we are interested in determining control
strategies that destroy ergodicity, or in other words the existence of an
attractive invariant measure.
Here, coupling arguments may be used as they provide criteria for the non-existence of a
unique invariant measure\protect\footnote{ Coupling arguments have been
  used since the theorem of Harris
  \citep{harris1960lower,lindvall2012lectures}, and are hence sometimes
  known as Harris-type theorems.  Generally, they link the existence of a
  coupling with the forgetfulness of initial conditions.}.


To formalise this discussion, let us denote the
space of trajectories of a $\bbSigma$-valued Markov chain $\{ X(k)\}_{k\in\N}$,
\emph{i.e.}, the space of all
sequences $(x(0),$ $x(1),$ $x(2), \ldots)$ with $x(k) \in \bbSigma$, ${k\in\N}$,
by $\bbSigma^\infty$ (the ``path space'').
Recall, for example, that $P_\lambda \in M(\bbSigma^\infty)$ is
the probability measure induced on
the path space by the initial distribution $\lambda$ of $X(0)$.
A coupling of two measures $P_{\mu_1},P_{\mu_2} \in M(\bbSigma^\infty)$ is a measure on $\bbSigma^\infty\times \bbSigma^\infty$ whose marginals coincide with $P_{\mu_1}, P_{\mu_2}$.
To be precise,
consider $\Gamma \in M(\bbSigma^\infty\times \bbSigma^\infty)$,
\emph{i.e.}, a measure over the product of space.
The measure $\Gamma$ can be projected to a measure over one or
the other factor space $\bbSigma^\infty$; we denote the projections by $\Phi^{(1)} \Gamma$ and $\Phi^{(2)} \Gamma$.
The set ${\mathbb C}(P_{\mu_1}, P_{\mu_2})$  of couplings
of $P_{\mu_1}, P_{\mu_2} \in M(\bbSigma^\infty)$ is then defined by
\begin{align*}
   \left \{
\Gamma \in M(\bbSigma^\infty \times \bbSigma^\infty) \; : \;
\Phi^{(1)} \Gamma = P_{\mu_1}, \Phi^{(2)} \Gamma  = P_{\mu_2}
\right \}.
\end{align*}
We say that a coupling $\Gamma \in {\mathbb{C}}$ is an asymptotic coupling if $\Gamma$ has
full measure on the pairs of convergent sequences. To make this precise
consider the set ${\mathbb H}$ defined by:
\begin{align*}
    \left \{ (x_1, x_2) \in \bbSigma^\infty \times \bbSigma^\infty \; : \; \lim_{k\to \infty} \left\| x_1(k) - x_2(k) \right\| =0
    \right \}.
\end{align*}
$\Gamma$ is an asymptotic coupling if $\Gamma({\mathbb H}) = 1$.
The following statement is a specialization of \citep[Theorem 1.1]{hairer2011asymptotic}
to our situation:

\begin{thm}[\cite{hairer2011asymptotic}]
\label{coupling-argument}
Consider an IFS with associated Markov operator $P$ admitting two ergodic invariant measures $\mu_1$ and $\mu_2$. The following are equivalent:
\begin{enumerate}
\item[(i)] $\mu_1 = \mu_2$.
\item[(ii)] There exists an asymptotic coupling of $P_{\mu_1}$ and $P_{\mu_2}$.
\end{enumerate}
\end{thm}
Consequently, if no asymptotic coupling of $P_{\mu_1}$ and $P_{\mu_2}$ exists,
then $\mu_1$ and $\mu_2$ are distinct.

\section{Main Results: Linear Systems}
\label{sec:comments-pi}

In this section, we present our main results for linear feedback systems, that is, systems in which agents, filter and controller are all linear systems. We begin by presenting some negative results for PI controllers; in fact, this first discussion is important for any control structure in which some of the components are marginally stable. We then focus on presenting conditions that ensure ergodicity to the closed-loop system -- our positive results.

\subsection{Negative Results: Controllers with Poles on the Unit Circle}
\label{sec:comments-pi-negative}

To illustrate the importance of the discussion on ergodicity we now present our first main result.
In many applications, controllers with integral action, such as the Proportional-Integral (PI) controller, are widely adopted \citep{Franklin,Franklin_dig}.
A simple PI control can be implemented as:
\begin{equation}\label{pid1}
\pi(k) = \pi(k-1) + \kappa \big[ e(k) - \alpha e(k-1) \big],
\end{equation}
which means its transfer function from $e$ to $\pi$, in terms of the ${\mathcal Z}$ transform, is given by
\begin{equation}\label{pid2}
C(z) = \kappa\frac{1 - \alpha z^{-1}}{1 - z^{-1}}.
\end{equation}
Since this transfer function is not asymptotically stable, any associated realisation matrix will not be Schur.
Note that this is the case for any controller with any sort of integral
action, \emph{i.e.}, a pole at $z = 1$.

\begin{thm}
    \label{thm:pole}
    Consider $N$ agents with states $x_i, i=1,\ldots,N$. Assume
    that there is an upper bound $m$ on the different values the agents
    can attain, \emph{i.e.}, for each $i$ we have $x_i \in \mathbb{A}_i =\{
    a_1,\ldots,a_{m_i} \}\subset \R$
    for a given set ${\mathbb{A}}_i$ and $1 \leq m_i \leq m$.

    Consider the feedback system in Figure \ref{system}, where ${\mathcal
      F}\, : \, y \mapsto \hat y$ is a finite-memory moving-average (FIR)
    filter.  Assume the controller ${\mathcal C}$ is a linear marginally
    stable single-input single-output (SISO) system with a pole $s_1 =
    e^{q j \theta}$ on the unit circle where $q$ is a rational number,
	$j$ is the imaginary unit,
	and $\theta$ is Archimedes' constant\protect\footnote{Here use $\theta$ for the Archimedes' constant of approximately $3.1416$ in order to avoid confusion with the feedback signal $\pi(k)$.}. In
    addition, let the
    probability functions $p_{il} : \R \to [0,1]$ be continuous for all
    $i=1, \ldots, N, l=1,\ldots,m_i$, \emph{i.e.}, if
    $\pi(k)$ is the output of ${\mathcal C}$ at time $k$, then $\probability(x_i(k+1)=a_l)=p_{il}(\pi(k))$. Then
    the following holds.

    \begin{enumerate}
      \item[(i)] The set ${\mathbb O}_{\mathcal F}$ of possible output values of the filter ${\mathcal F}$ is finite.
      \item[(ii)] If the real additive group ${\mathcal E}$ generated by
        $\{ r - \hat y \mid \hat y \in {\mathbb O}_{\mathcal F} \}$ with $r$ from \eqref{eq:1} is discrete, then
        the closed-loop system cannot be uniquely ergodic.
    \end{enumerate}
\end{thm}

\begin{rem}
One implication of the theorem it that
it is perfectly possible
for the closed loop both to perform its regulation function well and, at the same time, to
destroy the ergodic properties of the closed loop.
See Remark~\ref{rem:baseline} for a discussion of further implications.
Thereby, the performance of the closed loop needs to be studied both in terms of the classical regulation and in terms of the ergodic behaviour.
\end{rem}

\begin{pf}
(i)    By assumption, the states of the agents can only attain 
    finitely many values. Consequently, the set of possible values of $y$
    is finite, and thus also the set of possible outputs of the filter is
    finite, as it is just the moving average over a history of finite
    length.

(ii) We denote by ${\mathcal E}$ the additive subgroup of $\R$
    generated by the filter outputs.
By (i), the set of possible inputs to the linear part of the controller is
finite at any time $k\in \N$. Let $(A_c,B_c,C_c)$ be a minimal realization of
the linear controller with $A_c \in \R^{n_c \times n_c}$, $B_c,C_c^T\in \R^{n_c}$. Without  loss of generality, assume that
\begin{equation*}
    A_c =
    \begin{bmatrix}
        Q & 0 \\ 0 & R
    \end{bmatrix},\quad B_c =
    \begin{bmatrix}
        B_1 \\ B_2
    \end{bmatrix} , \quad C_c =
    \begin{bmatrix}
        C_1 & C_2
    \end{bmatrix}.
\end{equation*}
Here $Q$ is equal to $1,-1$ or a $2 \times 2$ orthogonal matrix with the
eigenvalues $s_1$ and $\overline{s_1}$. The matrix $R$ is marginally Schur
stable. We will concentrate on the first element(s) $\textrm{index}(x_c(k), 1)$ of the state $x_c$ of
the controller at time $k$ compatible with $Q$. That is: $\textrm{index}(x_c(k), 1) \in \mathbb{R}^1$ when $Q$ is a scalar,
and $\textrm{index}(x_c(k), 1) \in \mathbb{R}^2$ when $Q$ is a $2 \times 2$ matrix. Given an initial value
$x_c(0)$ and its first component $\textrm{index}(x_c(0), 1)$:
\begin{equation*}
    \textrm{index}(x_c(k), 1) = Q^k \textrm{index}(x_c(0), 1) + \sum_{\nu=0}^{k-1}Q^{k-\nu-1} B_1 e(\nu),
\end{equation*}
where the sequence $e(0), e(1), \ldots$ represents the input to the controller.
For some power $K\geq 1$ we have that $Q^K = I_2$ (or $1$), by assumption. We may thus
rearrange the sums and just consider finitely many powers of $Q$. This
induces a further summation over a subsequence of $\{ e(\nu) \}$, which by
construction lies in ${\mathcal E}$.
Thus
$\textrm{index}(x_c(k), 1)$ is an element of the set ${\mathbb Z}(x_{c}(0))$ given by
\begin{align*}
{\mathbb Z}(x_{c}(0)) := \Big \{ Q^k \textrm{index}(x_c(0), 1) +  \sum_{\nu=0}^{K-1} Q^\nu B_1 e(\nu) \\
\Big\vert \ k = 0,..,K-1 ,e(\nu) \in {\mathcal E} \Big \}.
\end{align*}
By assumption, this set is discrete in $\R$ or $\R^2$, as the case may
be. The state space of the controller may thus be partitioned into the
uncountably many equivalence classes under the equivalence relation on
$\R^{n_c}$ given by $x \sim y$, if the first element (for scalar $Q$ or first two element for $2 \times 2$ matrix $Q$) of $y$ is in ${\mathbb Z}(x)$,
i.e., $\textrm{index}(y, 1) \in {\mathbb Z}(x)$. These are
invariant under the evolution of the Markov chain. By
Proposition~\ref{prop:ergodic-car}~E2, ergodic invariant
measures are concentrated on one of these invariant sets.
Ergodic invariant measures that are concentrated on different equivalence
classes cannot couple asymptotically, as the respective
trajectories remain a positive distance apart. By Theorem~\ref{coupling-argument}, the
Markov chain cannot be uniquely ergodic. In particular, should there be
only one ergodic invariant measure $\mu$, then \eqref{eq:ergodicprop} cannot
hold for all deterministic initial conditions (just take a nonzero
continuous function $g$ which is zero on the support of $\mu$).
\end{pf}

While the conditions of the previous result appear fairly abstract, we would
like to point out that they apply in many practical settings. In an
implementation using standard digital computers all constants appearing
in the system description are rational numbers. It is therefore of
interest to observe that in this case the above theorem applies, as we note
in the next result.
\begin{cor}
    \label{cor:pole}
    In the situation of Theorem~\ref{thm:pole}, assume that ${\mathbb A}_i
    \subset \mathbb{Q}$ for all $i=1,\ldots,N$. Assume furthermore that $r
    \in \mathbb{Q}$ and that the coefficients of the FIR filter ${\mathcal F}$ are
    rational. Then the group ${\mathcal E}$ is discrete. If the  linear
    controller satisfies the assumptions of Theorem~\ref{thm:pole}, the closed-loop
    system cannot be uniquely ergodic.
\end{cor}
\begin{pf}
As $\mathbb{Q}$ is a field, it is easy to see that the set ${\mathbb O}_{\mathcal F}$ is contained in
$\mathbb{Q}$. Indeed, the possible outputs are obtained by manipulation of
rational numbers using linear maps with rational coefficients. It follows
that the finite set of generators of the additive group ${\mathcal E}$ is rational.
It follows that ${\mathcal E}$ is discrete and the final claim follows from
Theorem~\ref{thm:pole}\,(ii).
\end{pf}
\begin{rem}
We note that in real implementations it may happen that the common denominator
of the elements of ${\mathcal E}$ is so small that it is below machine
precision, which may lead to effects not predicted by
Corollary~\ref{cor:pole}.  But we do not pursue this question here.
\end{rem}

\begin{rem}
The inability to use integral action in managing the aggregate effect of populations clearly has profound implications for control design; in particular, for
rejecting certain types of disturbances.
\end{rem}

\begin{rem}
Note that control systems with integral action are already widely, and perhaps naively,
used in smart city applications,
including some of the present authors. While we shall not identify any specific miscreants (other than papers we have written ourselves) we note that, for example, in \citep{Arieh13}, precisely our setup is used to regulate the fleet emissions of a group of vehicles. Here, for a fleet of hybrid vehicles, $x_i(k)$ denotes whether or not vehicle $i$ switched into electric mode at time $k$ or not. Although measuring the pollution associated to each agent's action $x_i(k)$ is difficult, measurements $\hat y(k)$ of the aggregate level of pollution in a city are widely available from sensors. Then, based on this measurement, a central agency broadcasts a price signal $\pi(k)$, adjusted via a PI control, for example, based on the difference between a target level of pollution, and a measurement,  which is then used by the agents in order to probabilistically decide switching to electric mode or not. Problems of a similar nature arise in the context of balancing energy demand and supply in cities; see for example  \citep{ferc17}  \citep{192993}  \citep{Alagoz2013,Salsbury2013}. It is worth noting that the use of PI or PID control for such applications is entirely reasonable, with any issues around loss of ergodicity somewhat surprising and, from our perspective, unexpected.
\end{rem}

Before proceeding we now give a simple example to illustrate that the previous result
is not just of academic interest, but rather is also of  some practical importance.


{\bf Example: } {\em Let us illustrate the undesirable behaviour that may arise whenever a PI controller is being used in the closed-loop system.
In this example, we point out that the integral action may be heavily
dependent on the controller's initial state. In the
following we choose all data so that Corollary~\ref{cor:pole} is
applicable. That is, the sets ${\mathbb A}_i$, $r$ and the coefficients of
${\mathcal F}$ are chosen to be rational, and the probability functions $p_{ij}$
are continuous.

Consider the feedback system
depicted in Figure \ref{system} with $N = 10$ agents, whose states $x_i$ are in the set $\{0,1\}$; as before, if $x_i = 1$,
we say that agent $i$ has taken the resource or is {\em active}. In this example, we assume that only the first five agents start with the resource ($x_i(0) = 1$, $i = 1,\ldots,5$; $x_i(0) = 0$, $i = 6,\ldots,10$).

Our main goal is to regulate the number of active agents around the
reference value $r = 5$. We assume that five agents, namely $x_1$ to
$x_5$, have the following probabilities of being active ($i=1,\ldots,5$)
\begin{equation}
p_{i1}(x_i(k+1) = 1) = 0.02 +  \frac{0.95}{ 1 + \exp(-100(\pi(k) - 5))},
\notag
\end{equation}
whereas the remaining agents' probability of consuming the resource is
given by (for $i=6, \ldots, 10$)
\begin{equation}
p_{i1}(x_i(k+1) = 1) = 0.98 -  \frac{0.95}{ 1 + \exp(-100(\pi(k) - 1))}.
\notag
\end{equation}
As all agents have two options, we always have $p_{i0}=1-p_{i1}$.
Indeed, if the control signal $\pi(k) \gg 5$, then the first five agents are more likely to be active.
On the other hand, if $\pi(k) \ll 1$, then remaining ones are more likely to take the resource.

In this case, we implement two types of linear controllers
${\mathcal C}$: a PI controller and its lag approximant. The PI controller implements (\ref{pid1}) with $\kappa = 0.1$ and $\alpha = -4$.
This controller is approximated by a lag controller, whose transfer function is given by
\[
C(z) = \kappa \frac{1 - \alpha z^{-1}}{1 - \beta z^{-1}},
\]
with $\kappa = 0.1$, $\alpha = -4.01$ and $\beta = 0.99$.
The filter ${\mathcal F}$ is the moving average (FIR) filter defined by
\begin{equation}
  \hat y(k) = \frac{y(k) + y(k-1)}{2}.
\end{equation}

Our first observation from one simulation is that the filter output, $\hat y$, assumes, indeed,
a finite set of rational values, as shown in Figure \ref{sim_uma_so}. Hence, the conditions of Corollary~\ref{cor:pole}
are met by both the controller and the filter in the PI case. As the
closed-loop system is not uniquely ergodic in this case,
it is possible that undesirable characteristics may be observed during simulations; such behaviour should not be observed for the lag controller.

We simulated the feedback system depicted in Figure \ref{system} with the setup described above
and 2000 sample paths. The results of these Monte Carlo simulations are presented in Figures \ref{sim}--\ref{sim4}.
The figures use solid curves for the mean and shading for the area one standard deviation above and below the mean.
Figure \ref{sim} points out that the PI controller regulates the average number of active agents
$\bar y$, whereas the lag controller presents a steady-state error (as expected). However, Figure \ref{sim2} shows different
average trajectories of one of the five first agents, $\bar x_1$, for
different initial conditions of the controller ${\mathcal C}$, namely
$x_c(0) = 50$ and $x_c(0) = -50$. As the figure points out, this agent's
behaviour is completely dependent on the initial value of $x_c$, when
${\mathcal C}$ is the PI controller. It is important to note, however,
that this undesirable behaviour vanishes over the long run when a lag
controller is used.
That is, the system becomes uniquely ergodic and, hence, predictable.
We consider this unexpected dependency further in Figure \ref{sim4},
which points out the influence of the initial PI controller's state on
the average state of one of the first agents, including $\bar x_1$ over
the long run. Figure \ref{sim3} illustrates the dynamic response of the broadcast signal $\pi$ for both initial conditions
and both controllers; both cases converge to the same value for the lag structure and this is not observed
when PI is used.

\begin{figure}
    \begin{center}
   \includegraphics[width=\columnwidth]{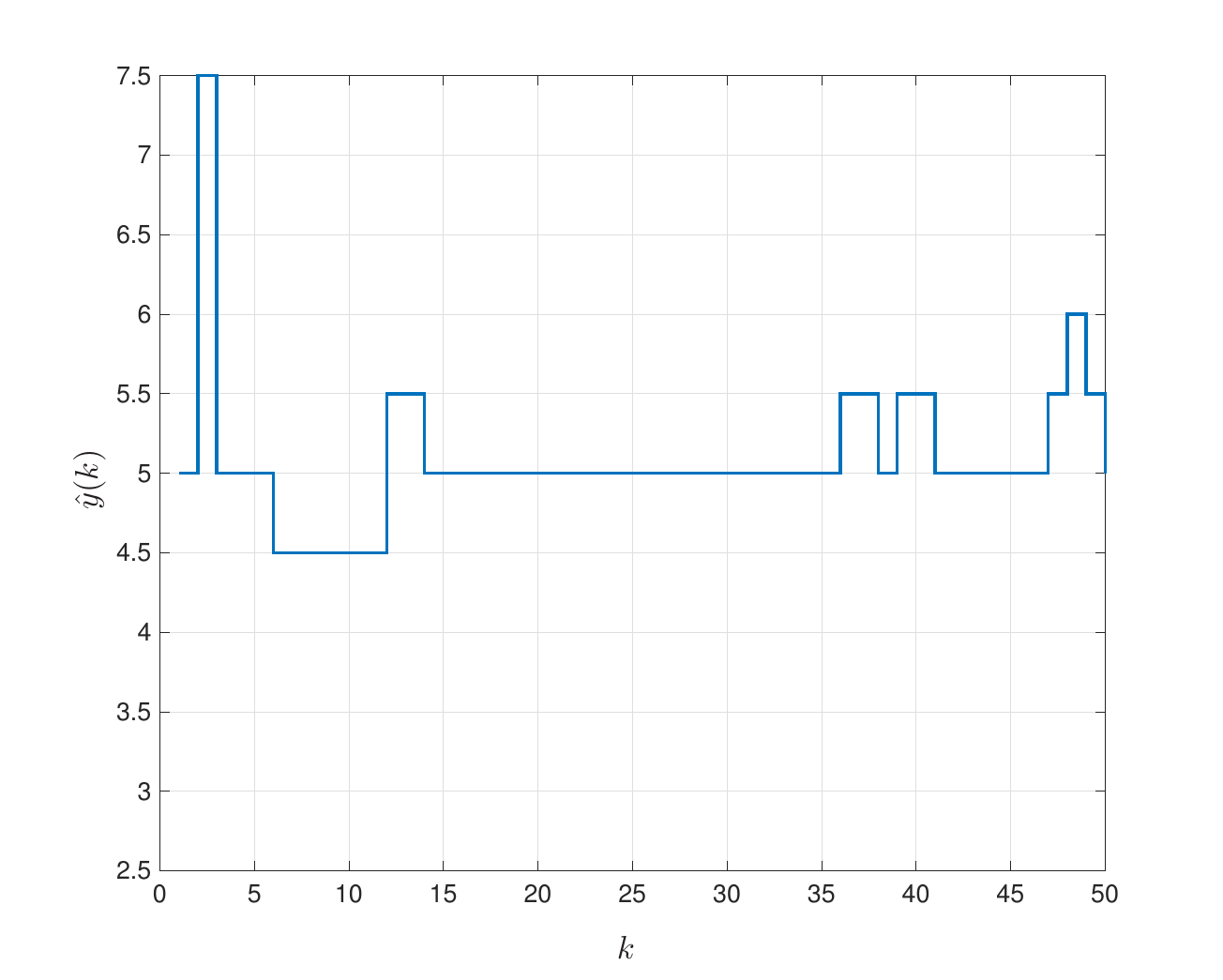}
    \caption{Filter output for a single simulation. Note that $\hat y(k)$ assumes a finite set of rational values,
		and hence satisfies conditions of Corollary~\ref{cor:pole}.}\label{sim_uma_so}
    \end{center}
\end{figure}

\begin{figure}
    \begin{center}
    \includegraphics[width=\columnwidth]{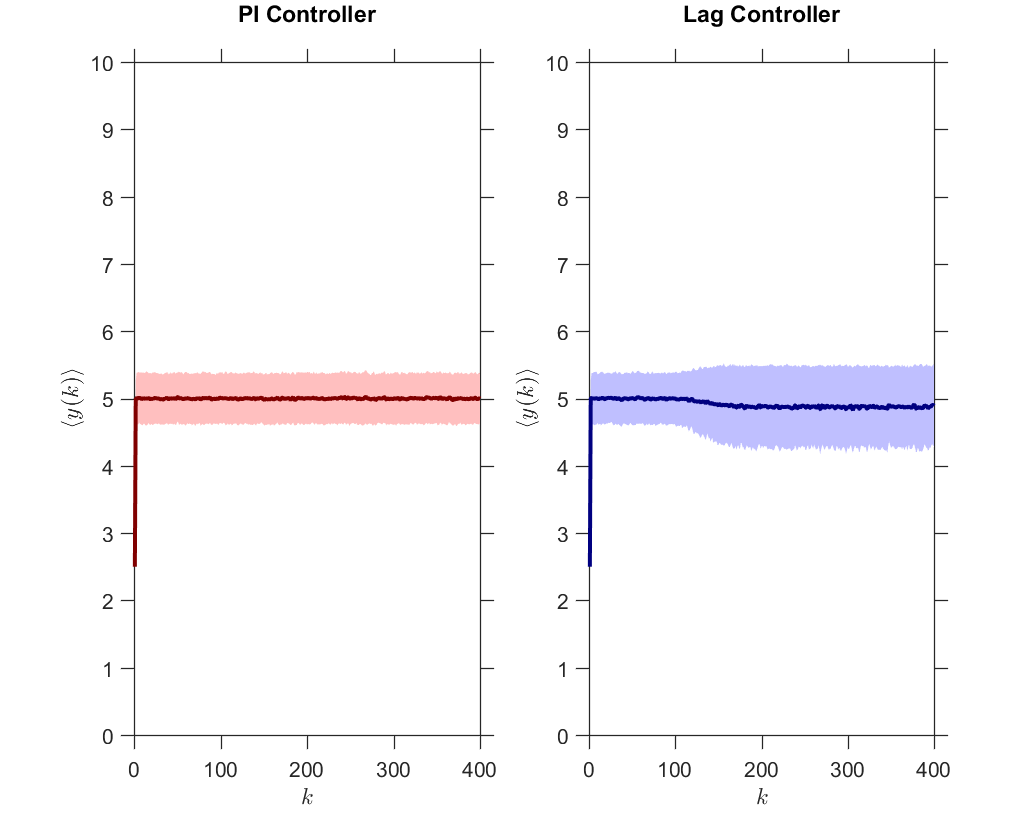}
    \caption{The solid curves represent the mean $y(k)$ over 2000 Monte Carlo simulations.
		The shaded areas denote one standard deviation around the mean.
		Regulation is observed (on average) for both the PI and for the lag controllers, within a given precision and for $x_c(0) = 50$.}\label{sim}
    \end{center}
\end{figure}

\begin{figure}
    \begin{center}
\includegraphics[width=\columnwidth]{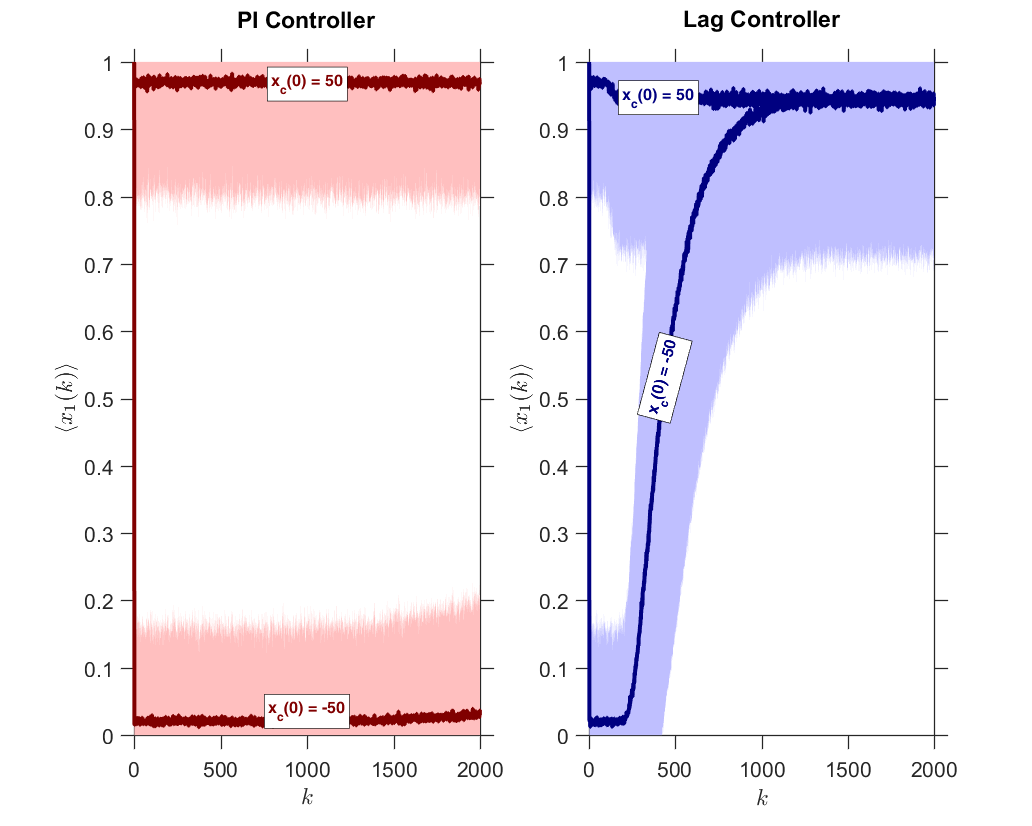}
    \caption{The solid curves represent the mean $x_1(k)$ over 2000 Monte Carlo simulations.
		The shaded areas denote one standard deviation around the mean.
		Four cases are considered in this figure, as we simulate the behaviour of $x_1$ for both controllers and two initial states, $x_c(0) = 50$ and $x_c(0) = -50$. Predictability is lost for the PI controller but is satisfied for the lag controller, as $x_c(0)$ does not affect $\langle x_1(k) \rangle$ for the second controller on the long run.}\label{sim2}
    \end{center}
\end{figure}

\begin{figure}
    \begin{center}
    \includegraphics[width=\columnwidth]{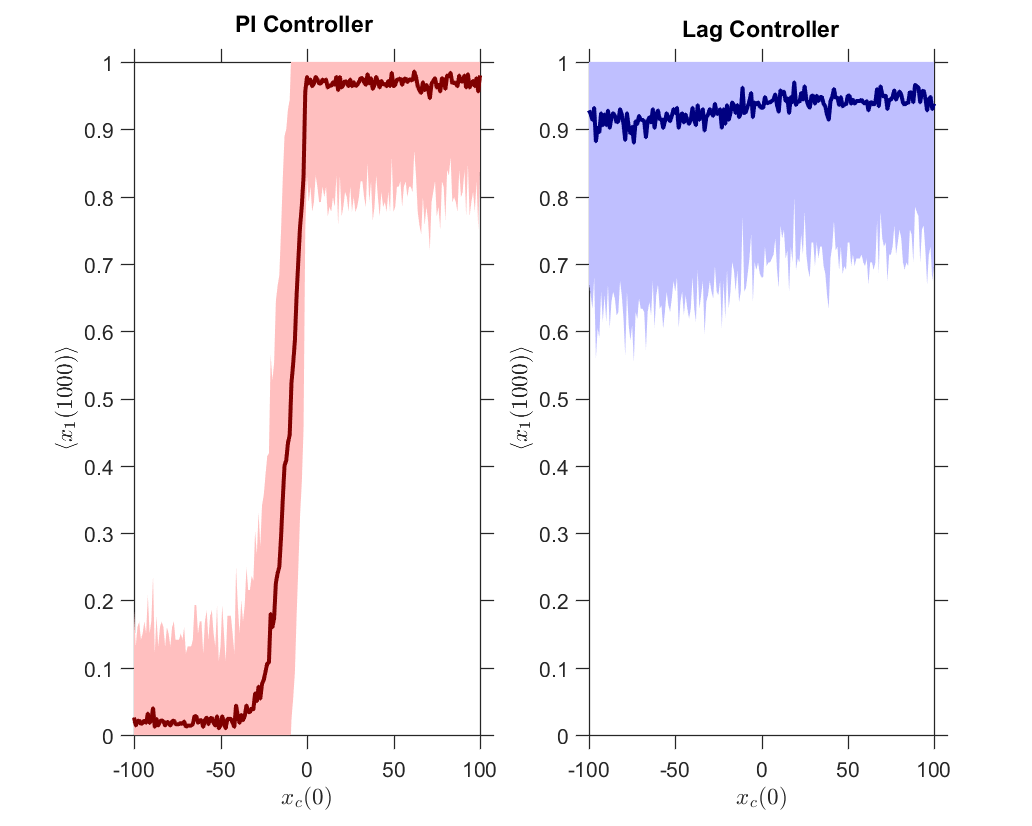}
    \caption{The solid curves represent the mean $x_1(1000)$  over 2000 Monte Carlo simulations.
				The shaded areas denote one standard deviation around the mean.
		This figure presents the behaviour of the average of
                $x_1(1000)$ for both controllers in dependence of the initial controller state
		$x_c(0)$. Predictability of $x_1(1000)$ is once again lost, when the PI controller is used.}\label{sim4}
    \end{center}
\end{figure}

\begin{figure}
    \begin{center}
    \includegraphics[width=\columnwidth]{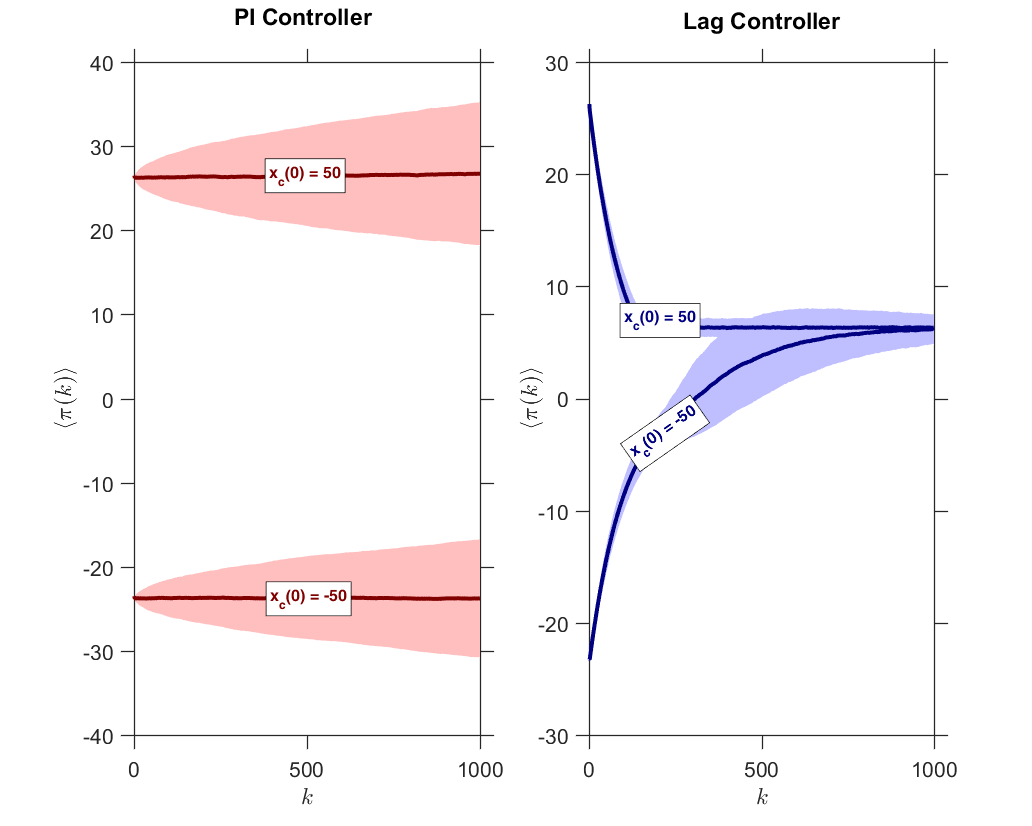}
    \caption{The solid curves represent the mean $\pi(k)$ over 2000 Monte Carlo simulations.
		The shaded areas denote one standard deviation around the mean.
		Four cases are considered in this figure: we consider both controllers for the simulations as well as two initial conditions,
		$x_c(0) = 50$ and $x_c(0) = -50$, for each of them.
		Once again, predictability is lost for the PI controller, as observed by changing $x_c(0)$.}\label{sim3}
    \end{center}
\end{figure}
}

\subsection{Positive Results: Ergodic Behaviour under Feedback}
\label{sec:comments-pi-positive}

Finally, to conclude this section, we note that coupling fails in our PI example
due to a lack of contractivity. Fortunately, for linear systems, the notion of
contractivity needed for unique ergodicity is relatively easy to enforce, and
we shall now provide conditions which guarantee a stable behavior,
for particular combinations of agent dynamics,
filter, and controller. Specifically, in the linear setting, the controller dynamics are:
\begin{equation} \label{eq_c}
{\mathcal C} ~:~ \left\{ \begin{array}{rcl}
x_c(k+1) & = & A_c x_c(k) + B_c e(k), \vspace{0.1cm} \\
\pi(k) & = & C_c x_c(k) + D_c e(k),
\end{array} \right.
\end{equation}
where $x_c \in\R^{n_c}$ is the internal state of the controller of dimension $n_c$.
We adopt a linear model for the $n_f$-dimensional filter ${\mathcal F}$, based on the classic IIR/FIR structures \citep{Oppenheim}.
Remembering that $y$ is the sum of each agents' output (as in Figure \ref{system}), one has
\begin{equation}\label{eq_f}
{\mathcal F} ~:~ \left\{ \begin{array}{rcl}
x_f(k+1) & = & A_f x_f(k) + B_f y(k),\vspace{0.1cm}\\
\hat y(k) & = & C_f x_f(k).
\end{array} \right.
\end{equation}
Finally, to be consistent with our discussion we consider populations of agents
 whose dynamics are described by:
\begin{equation}\label{eq:4:1}
{\mathcal S}_i ~:~ \left\{ \begin{array}{rcl}
x_i(k+1) & = & A_i x_i(k) + b_i, \vspace{0.1cm} \\
y_i(k) & = & c_i^T x_i(k) + d_i,
\end{array}\right.
\end{equation}
where $A_i \in \R^{n_i \times n_i}$, $c_i \in \R^{n_i}$. The inputs $b_i$ and $d_i$ are random variables that take values in $\R^{n_i}$ and $\R$,
respectively, with $\mathbb{P}(b_i = b_{ij}) = p_{ij}(\pi)$ and $\mathbb{P}(d_i = {d}_{ij}) = p_{ij}'(\pi)$. Note that this is a generalisation
of the situation when agents switch between two states {on, off}.
The following result gives conditions for unique ergodicity.

\begin{thm} \label{thm01}
Consider the feedback system depicted in Figure \ref{system}, with ${\mathcal C}$ and ${\mathcal F}$ given in (\ref{eq_c}) and (\ref{eq_f}). 
Assume that each agent  
 $i \in \{1,\cdots,N\}$ has state $x_i$ with dynamics governed by the
 affine stochastic difference equations given in \eqref{eq:4:1}, where
 $A_i$ are Schur matrices and $b_i$ and $d_i$ are chosen, at each time
 step, from the sets $\{b_{ij}\} \subset \R^{n_i}$ and $\{d_{i\ell}\}
 \subset \R$ according to Dini continuous probability functions $p_{ij}(\cdot)$, respectively $p'_{i\ell}(\cdot)$, that verify \eqref{eq:problaws}.
 Assume furthermore that there are scalars $\delta, \delta' > 0$ such that
$p_{ij}(\pi) \geq  \delta > 0$, $p'_{ij}(\pi) \geq  \delta' > 0$ for all
$(i,j)$ and all $\pi \in \Pi$. Then, for every stable linear controller
$\mathcal{C}$ and every stable linear filter $\mathcal{F}$ compatible with
the system structure, the feedback loop converges in distribution to a unique invariant measure.
\end{thm}

\begin{pf}
Following \citep{BarnsleyDemkoEltonEtAl1988}, the proof is centred at the construction of an iterated function system (IFS) with
place-(state-)dependent probabilities that describes the feedback system. To this end, consider the augmented state
$\xi \defeq [x^T~,~x_f^T~,~x_c^T]^T \in \mathbb{X}$, whose dynamic behaviour is described by the difference equation
\begin{equation}
\label{dynamical}
\xi(k+1) = {\mathcal W}_\ell(x) \defeq \mathcal{A}\xi(k) + \beta_\ell,
\end{equation}
where
\begin{equation}
\label{mathcaldef}
\mathcal{A} \defeq
\begin{bmatrix}
  \hat A & 0 & 0 \\
  B_f {\bf 1}^T \hat C & A_f & 0\\
  0 & -B_c C_f & A_c
\end{bmatrix}
\end{equation}
\normalsize
where $\mathbf{1}$ is the vector of ones, $\hat A \defeq \mathbf{diag}(A_i)$, $\hat C \defeq
\mathbf{diag}(c_i^T)$ and $\beta_\ell$ is built from all the combinations of the
vectors $b_{ij}$, the scalars ${d}_{ij}$ and other signals. To apply Corollary 2.3 from
\citep{BarnsleyDemkoEltonEtAl1988}, two observations must be made. First,
note that each map ${\mathcal W}_\ell$ is chosen with probability $p_\ell(\pi) \geq
\prod_{i=1}^N \delta_i > 0$ and thus these probabilities are bounded away from
zero. Second, since $\sigma({\mathcal A}) = \sigma(\hat A) \cup
\sigma(A_f) \cup \sigma(A_c)$ and, by hypothesis, $A_i$, $A_f$
and $A_c$ are Schur matrices,  for any induced matrix norm there exists $m
\in \N$ sufficiently large such that $\|{\mathcal A}^m\| < 1$. This
provides the required average contractivity after $m$ steps. The result then follows from \citep{BarnsleyDemkoEltonEtAl1988}. The proof is complete.
\end{pf}

\begin{rem}
Dini's condition on the probabilities may be replaced by simpler, more conservative assumptions, such as Lipschitz or H\"older conditions \citep{BarnsleyDemkoEltonEtAl1988}. Also the requirement $p_{ij}(\pi) \geq \delta_i > 0$ in the theorem statement is not an artefact of our analysis, as $p_{ij}(\pi) = 0$ may lead to a non-ergodic behaviour.
\end{rem}

\begin{rem}
As we have stressed in Remark~\ref{rem:baseline}, the existence of an attractive invariant measure is only a baseline ergodic property.
Under mild but technical additional assumptions, it is possible to prove a geometric rate of convergence
\citep{steinsaltz1999}.
Under further assumptions, one can prove results concerning the ``shape'' of the unique invariant measure,
where it exists, such as
moment bounds  \citep{walkden2007invariance}.
These would complement the present theorem.
\end{rem}

\begin{rem}
In Theorem \ref{thm01}, we assume that the feedback loop only consists of stable systems.
This assumption may seem quite restrictive, but as shown in the examples, even marginal stability in
one of its parts may lead to completely unpredictable results. As the
result from \citep{BarnsleyDemkoEltonEtAl1988} only requires contractivity
on average, more general results are possible, but beyond the scope of the
present paper. We also point out that each agent may have a local stabilising controller,
as we are analysing the feedback loop from a macroscopic level.
\end{rem}

\section{Extensions: Non-linear systems and discrete action spaces}
\label{sec:switched-linear-non}

Considering that Theorem \ref{Barnsley} does not require linearity,
it is clear that one can extend its use to non-linear systems under suitable assumptions.

\subsection{Non-linear controllers}
\label{sec:nonlinear}

A particular case of the general setup described in  Figure \ref{system}
is given by systems of the following form:
\begin{align}
    \label{eq:nonlinear-agents}
\phantom{{\mathcal A}}  &\left\{ \begin{array}{ccl}
    x_i(k+1) &\in& \{  {\mathcal W}_{ij}(x_i(k)) \; \vert \; j=1,\ldots, w_i \}  \\
    y_i(k) &\in& \{  {\mathcal H}_{ij}(x_i(k))  \; \vert \; j=1,\ldots, h_i \}, \\
\end{array} \right. \\
    & \hspace*{1.6cm} y(k) = \sum_{i=1}^N y_i(k),
\end{align}
\begin{equation}
    \label{eq:nonlinear-filter}
{\mathcal F} ~:~ \left\{ \begin{array}{ccl}
 x_f(k+1) &=& {\mathcal W}_{f}(x_f(k),y(k)) \\
    \hat y(k) &=& {\mathcal H}_{f}(x_f(k),y(k)),
\end{array} \right.
\end{equation}
\begin{equation}
    \label{eq:nonlinear-cont}
 {\mathcal C} ~:~ \left\{ \begin{array}{ccl}
   x_c(k+1) &=& {\mathcal W}_{c}(x_c(k),\hat y(k),r) \\
    \pi(k) &=& {\mathcal H}_{c}(x_c(k),\hat y(k),r),
\end{array} \right.
\end{equation}
In addition, we have Dini continuous probability functions
$p_{ij},p'_{il}:\Pi \to [0,1]$ so that
the probabilistic laws \eqref{eq:problaws} are satisfied.
%
If we denote by $\mathbb{X}_i, i=1,\ldots,N, \mathbb{X}_C$ and $\mathbb{X}_F$ the state spaces of the
agents, the controller and the filter, then the system evolves on the
overall state space $\mathbb{X} := \prod_{i=1}^N \mathbb{X}_i \times \mathbb{X}_C \times \mathbb{X}_F$
according to the dynamics
\begin{equation}
    \label{eq:6}
    x(k+1) := \begin{pmatrix}
        (x_i)_{i=1}^N  \\
        x_f \\
        x_c
    \end{pmatrix} (k+1) \in \{ F_m(x(k)) \,\vert\, m \in {\mathbb M}\}.
\end{equation}
where each of the maps $F_m$ is of the form
\begin{equation}
    \label{eq:77}
F_m(x(k)) \defeq  \begin{pmatrix}
        ( {\mathcal W}_{ij}(x_i (k)) )_{i=1}^N  \\
        {\mathcal W}_f(x_f(k), \sum_{i=1}^N  {\mathcal H}_{i\ell}(x_i (k))) \\
        {\mathcal W}_c(x_c(k), {\mathcal H}_f(x_f(k), \sum_{i=1}^N  {\mathcal H}_{i\ell}(x_i (k))))
    \end{pmatrix}
\end{equation}
and the maps $F_m$ are indexed by indices $m$ from the set
\begin{equation}
    \label{eq:8}
    \mathbb{M} \defeq \prod_{i=1}^N \{ (i,1), \ldots, (i,w_i) \} \times \prod_{i=1}^N \{ (i,1), \ldots, (i,h_i) \}.
\end{equation}
By the independence assumption on the choice of the transition maps and
output maps for the agents, for each multi-index
$m=((1,j_1),\ldots,(N,j_N),(1,l_1),\ldots,(N,l_N))$ in this set, the
probability of choosing the corresponding map $F_m$ is given by
\begin{multline}
    \label{eq:9}
    \probability \left( x(k+1)=F_m(x(k)) \right) = \\ \left(\prod_{i=1}^N
      p_{ij_i}(\pi(k)) \right) \left( \prod_{i=1}^N p'_{il_i}(\pi(k))
    \right) =: q_m(\pi(k)).
\end{multline}

\begin{thm} \label{thm02}
Consider the feedback system depicted in Figure \ref{system}.
Assume that each agent
 $i \in \{1,\cdots,N\}$ has a state governed by the non-linear iterated
 function system
\begin{align}
    \label{eq:nonlinear}
    x_i(k+1) &= {\mathcal W}_{ij}(x_i(k)) \\
    y_i(k) &= {\mathcal H}_{ij}(x_i(k)),
\end{align}
where ${\mathcal W}_{ij}$ and ${\mathcal H}_{ij}$ are
globally Lipschitz-continuous functions with
global Lipschitz constant $l_{ij}$, resp. $l'_{ij}$.
Assume we have Dini continuous probability functions
$p_{ij},p'_{il}$ so that
the probabilistic laws \eqref{eq:problaws} are satisfied. Assume
furthermore that there are  scalars $\delta, \delta' > 0$ such that
$p_{ij}(\pi) \geq  \delta > 0$,
$p'_{ij}(\pi) \geq  \delta' > 0$ for all $(i,j)$.
 Further, assume that the following contractivity condition holds:
for all $1 \le i \le N, 1 \le j \le J$: $l_{ij} < 1$.
Then, for every stable linear controller $\mathcal{C}$ and every stable
linear filter $\mathcal{F}$ compatible with the feedback structure, the feedback loop has a unique attractive
invariant measure. In particular, the system is uniquely ergodic.
\end{thm}

\begin{pf}
Similarly to the proof of Theorem~\ref{thm01} the assumptions on the
Lipschitz constants and the the internal asymptotic stability of
controller and filter guarantees that suitable iterates of the maps $F_m$
are strict contractions. The result then follows again from Theorem 2.1 and Corollary 2.2 of
 \citep{BarnsleyDemkoEltonEtAl1988}.
\end{pf}

\begin{rem}
\label{rem:QUAD}
Notice that Lipschitz continuity can be rephrased in many ways.
For instance, there is the QUAD condition \citep{5641620}, the sector condition, or, when restricting to convex functions,  the bounded subgradient condition.
The sector condition requires that there exist constants $c_1$ and $c_2$ such that the vector-valued functions
${\mathcal W}(x) := [{\mathcal W}_i (x)]$ and ${\mathcal H}(x):= [{\mathcal H}_i (x)]$ satisfy
${\mathcal W}(x)^T[{\mathcal W}(x) - c_1 x] \leq 0$
and
${\mathcal H}(x)^T[{\mathcal H}(x) - c_2 x] \leq 0$.
The bounded subgradient condition requires the existence of constants
$c_3, c_4$, such that for a given norm $|\cdot|$, for all $z, z'$ in the subdifferentials of ${\mathcal W}, {\mathcal H}$, respectively,
at all points in the domains of the respective functions,
we have that $|z| \le c_3, |z'| \le c_4$.
Here the functions ${\mathcal W}(x)$, ${\mathcal H}(x)$ are assumed to be convex with a non-empty subdifferential throughout their domains,
but not necessarily differentiable.
The equivalence follows from basic convex analysis, e.g., as a corollary of Lemma 2.6 in \citep{shalev2012online}.
\end{rem}

Let us refer back to Remark~\ref{rem:baseline} for a discussion of the
importance of unique ergodicity;
notably, the existence of the limit \eqref{eq:ergodicprop} almost surely.

\subsection{Discrete Action Spaces}
\label{sec:discrete}

Next, consider the case, when the agents' actions are
limited to a finite set. In this case the Lipschitz conditions in
Theorem~\ref{thm02} cannot be satisfied except in trivial cases. In this
case we may use results in \citep{werner2005contractive,werner2004ergodic}
to obtain ergodicity results.

The general setup of the following result is that for each agent $i$ the
set ${\mathbb A}_i$ is finite. Then, $\mathbb{X}_S := \prod_{i=1}^N {\mathbb A}_i$
is finite and we consider the directed graph $G=(\mathbb{X}_S, E)$, where there is an
arc between vertices representing $(x_i) \in \mathbb{X}_S$ and $(y_i) \in \mathbb{X}_S$, if
there is a choice of maps ${\mathcal W}_{ij}$ in \eqref{eq:nonlinear-agents} such that $({\mathcal W}_{ij}(x_i)) = (y_i)$.

\begin{thm}{}
    \label{thm03}
    Consider the feedback system depicted in Figure \ref{system}.  Assume
    that ${\mathbb A}_i$ is finite for each $i$. Assume that each agent $i
    \in \{1,\ldots,N\}$ has a state governed by the non-linear stochastic
    difference equations \eqref{eq:nonlinear}. Assume we have Dini
    continuous probability functions $p_{ij},p'_{il}$ so that the
    probabilistic laws \eqref{eq:problaws} are satisfied. Assume
    furthermore that there are scalars $\delta, \delta' > 0$ such that
    $p_{ij}(\pi) \geq \delta > 0$, $p'_{ij}(\pi) \geq \delta' > 0$ for all
    $(i,j)$ and all $\pi$. Then, for every stable linear controller $\mathcal{C}$ and
    every stable linear filter $\mathcal{F}$ the following holds:

If the graph
$G=(\mathbb{X}_S, E)$ is strongly connected, then there exists an
invariant measure for the feedback loop. If in addition, the adjacency matrix of
the graph is primitive, then the invariant measure is attractive and the
system is uniquely ergodic.
\end{thm}

\begin{pf}
    This is a consequence of \citep{werner2004ergodic} and the observation
    that the necessary contractivity properties follow from the internal
    asymptotic stability of controller and filter.
\end{pf}

We note that a simple condition for the primitivity of the graph
$G=(\mathbb{X}_S, E)$ is that for each agent the graph describing the
possible transitions is primitive.

\begin{rem}
    We note that there are a few cases, to which both Theorem~\ref{thm01}
    and Theorem~\ref{thm03} are applicable. Namely, if in the assumptions
    of Theorem~\ref{thm01} all $A_i=0$, then each agent at each time step
    chooses its next state $x_i(k+1)$ independently of the current state
    $x_i(k)$. In Theorem~\ref{thm03} this corresponds to the case that the
    ${\mathcal W}_{ij}$ are constant maps and the directed transition graph of the system is
    complete.
\end{rem}

\section{Conclusions and Further Work}

Within feedback systems, the control of ensembles of agents presents a particularly challenging area for further study.
Practically important examples of such systems arise in Smart Cities.
Typically, such problems deviate from classical control problems in two main ways.
First, even though ensembles are typically too large to allow for a microscopic approach,
they are not sufficiently large to allow for a meaningful fluid (mean-field) approximation.
Second, the regulation problem concerns not only the ensemble, but also the individual agents; a certain quality of service should be provided to each agent.
We have formulated this problem as an iterated function system with the objective of designing an ergodic control,
and demonstrated that controls with poles on the unit circle (e.g., PI) may destroy ergodicity even for benign ensembles.


\begin{ack}                               
This work was in part supported by Science Foundation Ireland grant 16/IA/4610, in part by
 Funda\c{c}\~ao de Amparo \`a Pesquisa do Estado de S\~ao Paulo (FAPESP) grants 2016/19504-7 and 2018/04905-1, in part by Conselho Nacional de Desenvolvimento Cient\'ifico e Tecnol\'ogico (CNPq) grant 305600/2017-6
 and in part by the European Union Horizon 2020 Programme (Horizon2020/2014-2020), under grant agreement no.\ 68838.
\end{ack}

\bibliography{ref,mdps}

\clearpage
\onecolumn
\appendix
\section{An Overview of Notation}
\label{app:symbols}

In general, our notation employs the following rules:
Upper-case letters are used for matrices, (in caligraphic) groups, maps, operators, (and in blackboard) sets, and spaces,
while lower-case letters are used for vectors, scalars, and functions.
Subscripts are used to distinguish symbols; time-indexed symbols are followed by the time index in parentheses, as in $x(k)$.
Superscript is used for exponentiation and, in $^T$, for the transpose operation.
Caligraphic fonts are used for maps and operators.
Blackboard fonts are used for sets, spaces, and the probability operator $\mathbb{P}$.
In particular, sets of real, rational, and natural numbers are indicated by $\mathbb{R}$, $\mathbb{Q}$, and $\mathbb{N}$, respectively.
For a given space $\mathbb{X}$, $M(\mathbb{X})$ indicates the set of all probability measures over $\mathbb{X}$ and $\mathbb{X}^\infty$
denotes its associated path space, which consists of infinite right-sided sequences over $\mathbb{X}$.

In the following table of notations, we list related groups of symbols in the order of appearance;
within each group of symbols, symbols are ordered alphabetically, Latin first and Greek second.

\begin{tabularx}{\linewidth}{ l | X }
 \caption{A Table of Notation}\\\toprule\endfirsthead
    \toprule\endhead
    \midrule\multicolumn{2}{l}{\itshape continues on next page}\\\midrule\endfoot
    \bottomrule\endlastfoot
    \textbf{Symbol} & \textbf{Meaning} \\\midrule
\multicolumn{2}{l}{Standard sets:} \\
	\midrule
$\mathbb{N}$ & the set of natural numbers \\
$\mathbb{Q}$ & the set of rational numbers \\
$\mathbb{R}$ & the set of real numbers  \\
\multicolumn{2}{l}{Constants:} \\
\midrule
$\mathbf{1}$ & a compatible vector of ones \\
$c_1,\cdots,c_4$ & constants used in Remark \ref{rem:QUAD}\\
$h_i$ & the number of output maps ${\mathcal H}_{ij}$ \\
$l_{ij}$ &  Lipschitz constant for a transition map \\
$l'_{ij}$ & Lipschitz constant for an output map \\
$N$ & the number of agents \\
$n$ & a dimension of a generic state space $\bbSigma$ \\
$n_c$ & dimension of the state of the controller \\
$n_f$ & dimension of the state of the filter \\
$n_i$ & dimension of the state of $i$th agent's private state  \\
$m_i$ & the number of possible actions of agent $i$ \\
$m$ & an upper bound on the number of possible actions of any agent \\
$r$ & the reference value, i.e., desired value of $y(k)$ \\
$\overline{r}_i$ & $i$th agent's expected share of the resource over the long run \\
$w_i$ & the number of state transition maps ${\mathcal W}_{ij}$ \\
$z$ & ${\mathcal Z}$-transform variable \\
$\alpha$ & a constant used in the PI controller \eqref{pid1} or its lag approximant  \\
$\beta$ & a constant used in a lag controller  \\
$\delta$ & a constant in the iterated function system ergodicity used in Theorem~\ref{Barnsley}  \\
$\delta'$ & a constant in the iterated function system ergodicity used in Theorem~\ref{thm01} \\
$\eta$ & a lower bound on the values of probability functions \\
$\kappa$ & a constant used in the PI controller \eqref{pid1} or its lag approximant  \\
\multicolumn{2}{l}{Blocks:} \\
\midrule
$\mathcal{C}$ & controller representing the central authority \\
${\mathcal F}$ & filter  \\
$\mathcal{S}_1$, $\mathcal{S}_2$, \ldots, $\mathcal{S}_N$ & systems modelling agents \\
\multicolumn{2}{l}{Maps and operators:} \\
\midrule
$f_j$ & a generic map in a generic iterated function system \\
$g$ & a generic function of \eqref{eq:ergodicprop} \\
${\mathcal H}_{ij}$ & an output map \\
${\mathcal H}_{c}$ & a map modelling the controller in \eqref{eq:nonlinear-cont} \\
${\mathcal H}_{f}$ & a map modelling the filter in \eqref{eq:nonlinear-filter} \\
$P(x,\mathbb{G})$ & a generic transition operator \\
$P$ & a state-and-signal-to-state transition operator \\
$\mathcal P_k$ & a measure-space-over-states-to-measure-space-over-states operator \\
$p_j$ & a probability function of a generic iterated function system \\
$p_{ij}$   & a probability function for the choice of agent $i$'s transition map \\
$p'_{i\ell}$  & a probability function for the choice of agent $i$'s output map \\
${\mathcal W}_{ij}$ & transition maps \\
${\mathcal W}_{c}$ & maps modelling the controller in \eqref{eq:nonlinear-cont} \\
${\mathcal W}_{f}$ & maps modelling the filter in \eqref{eq:nonlinear-filter} \\
$\Phi^{(1)}$ & projector from a measure over the product of the two path spaces to a single path space \\
$\Phi^{(2)}$ & projector from a measure over the product of the two path spaces to a single path space \\
\multicolumn{2}{l}{Deterministic matrices used in the definitions:} \\
\midrule
$A_i$  & a matrix used in the agent dynamics \eqref{eq:4:1} of Theorem~\ref{thm01} \\
$A_c$ & a matrix used in the controller \eqref{eq_c} of Theorem~\ref{thm01} \\
$A_f$ & a matrix used in the filter \eqref{eq_f} of Theorem~\ref{thm01}  \\
$B_c$ & a matrix used in the controller \eqref{eq_c} of Theorem~\ref{thm01} and in Theorem~\ref{thm:pole} \\
$B_f$ & a matrix used in the filter \eqref{eq_f} of Theorem~\ref{thm01}  \\
$B_{i}$ & a matrix used in the agent dynamics \eqref{eq:4:1} of Theorem~\ref{thm01} \\
$C_c $ & a matrix used in the controller \eqref{eq_c} of Theorem~\ref{thm01} and in Theorem~\ref{thm:pole} \\
$C_f$ & a matrix used in the filter \eqref{eq_c} of Theorem~\ref{thm01}   \\
$c_i$   & a vector used in the agent dynamics \eqref{eq:4:1} of Theorem~\ref{thm01} \\
$D_c$ & a matrix used in the controller \eqref{eq_c} of Theorem~\ref{thm01} \\
$d_i$   & a vector used in the agent dynamics \eqref{eq:4:1} of Theorem~\ref{thm01} \\
$X(k)$ & element of a generic state space  \\
$\{ X(k)\}_{k\in\N}$ & a generic Markov chain \\
\multicolumn{2}{l}{Deterministic matrices and constants used in proofs:} \\
\midrule
$\mathcal{A}$ & augmented state transition matrix in \eqref{mathcaldef} \\
$\hat A$ & $\mathbf{diag}(A_i)$ in \eqref{mathcaldef} \\
$\beta_\ell$ & is built from all $b_{ij}$, ${d}_{ij}$, and other signals in \eqref{mathcaldef} \\
$\hat C$ & $\mathbf{diag}(c_i^T)$ in \eqref{mathcaldef} \\
$E$ & edges of a graph in the proof of Theorem \ref{thm03} \\
$G$ & a graph in the proof of Theorem \ref{thm03} \\
$K$ & a constant used in the proof of Theorem~\ref{thm:pole} \\
$Q$ & a scalar or matrix used in the proof of Theorem~\ref{thm:pole} \\
$R$ & a marginally Schur matrix used in the proof of Theorem~\ref{thm:pole} \\
\multicolumn{2}{l}{Algebras, groups, sets, and spaces:} \\
\midrule
$\mathbb{A}_i$ & the set of $i$th agent's actions \\
${\mathbb C}$ & set of couplings  \\
$\mathbb{D}_i$ & set of possible resource demands of agent $i$ \\
$\mathbb{B}(\bbSigma)$ & Borel $\sigma$-algebra  \\
${\mathcal E}$ & a real additive group of Theorem~\ref{thm:pole} \\
$\mathbb{G}$ & a generic event \\
$\mathbb{H}$ & a set used in the definition of asymptotic couplings \\
$\mathbb{J}$ & index set of a generic iterated function system \\
${\mathbb M}$ & an index set of \eqref{eq:8}   \\
$M(\mathbb{X})$ & a measure-space over $\mathbb{X}$ \\
$M(\mathbb{X}^{\infty})$ & a measure space over the path space \\
${\mathbb O}_{\mathcal F}$ & set of possible output values of the filter ${\mathcal F}$  \\
$\bar{\mathbb{S}}$ & a measurable subset of a state space \\
$\mathbb{X}_i$ & a private state space of agent $i$, often $\R^{n_i}$ \\
$\mathbb{X}_S$ & a product space of state spaces of all agents \\
$\mathbb{X}_F$ & a space of internal states of the filter \\
$\mathbb{X}_C$ & a space of internal states of the central controller \\
$\mathbb{X}$ & a state space of the controller, the filter, and the agents, combined \\
$\mathbb{X}^{\infty}$ & a path space \\
${\mathbb Z}$ & set used in the proof of Theorem~\ref{thm:pole} \\
$\bbPi$ & the set of admissible broadcast control signals \\
$\bbSigma$ & a generic state space \\
\multicolumn{2}{l}{Random variables:} \\
\midrule
$e(k)$ & the error signal at time $k$, i.e., $\hat y(k) - r$ \\
$\pi(k)$ & the signal broadcast at time $k$ \\
$x_c(k)$ & the internal state of the controller at time $k$ of dimension $n_c$ \\
$x_f(k)$ & the internal state of the filter at time $k$ of dimension $n_f$ \\
$x_i(k)$ & the internal state of agent $i$ at time $k$ \\
$y_i(k)$ & the resource utilisation of agent $i$ at time $k$ \\
$y(k)$   & the aggregate resource utilisation at time $k$ \\
$\hat y(k)$ & the value of $y(k)$ filtered by filter $\mathcal{F}$ \\
$\lambda$ & initial state (distribution) \\
\multicolumn{2}{l}{Further notation} \\
\midrule
$P_\lambda$ & a probability measure induced on the path space  \\
$\mu$ & a generic measure, usually on state space $\bbSigma$  \\
$\Gamma$ & a generic measure over the product of the two path spaces, potentially a coupling \\
$\nu$ & a counter
\end{tabularx}

\end{document}